\documentclass[11pt]{article}
\usepackage{amsfonts}
\usepackage{amsthm}
\usepackage{amsmath}
\usepackage[all]{xypic}
\makeatletter
\def\newtheorem#1{\@ifnextchar[{\@othm{#1}}{\@nthm{#1}}}
\def\@nthm#1#2{%
\@ifnextchar[{\@xnthm{#1}{#2}}{\@ynthm{#1}{#2}}}
\def\@xnthm#1#2[#3]{\expandafter\@ifdefinable\csname #1\endcsname
{\@definecounter{#1}\@newctr{#1}[#3]%
\expandafter\xdef\csname the#1\endcsname{\expandafter\noexpand
  \csname the#3\endcsname \@thmcountersep \@thmcounter{#1}}%
\global\@namedef{#1}{%
  \@thm{#1}{#2}}\global\@namedef{end#1}{\@endtheorem}}}
\def\@ynthm#1#2{\expandafter\@ifdefinable\csname #1\endcsname
{\@definecounter{#1}%
\expandafter\xdef\csname the#1\endcsname{\@thmcounter{#1}}%
\global\@namedef{#1}{%
  \@thm{#1}{#2}}\global\@namedef{end#1}{\@endtheorem}}}
\def\@othm#1[#2]#3{%
  \@ifundefined{c@#2}{\@nocounterr{#2}}%
  {\expandafter\@ifdefinable\csname #1\endcsname
  {\global\@namedef{the#1}{\@nameuse{the#2}}%
\global\@namedef{#1}{\@thm{#2}{#3}}%
\global\@namedef{end#1}{\@endtheorem}}}}
\def\@thm#1#2{\refstepcounter
    {#1}\@ifnextchar[{\@ythm{#1}{#2}}{\@xthm{#1}{#2}}}
\def\@xthm#1#2{\@begintheorem{#2}{\csname the#1\endcsname}\ignorespaces}
\def\@ythm#1#2[#3]{\@opargbegintheorem{#2}{\csname
       the#1\endcsname}{#3}\ignorespaces}
\def\@thmcounter#1{\noexpand\arabic{#1}}
\def\@thmcountersep{.}
\def\@begintheorem#1#2{\trivlist
   \item[\hskip \labelsep{\bfseries #2\ #1.}]\itshape}
\def\@opargbegintheorem#1#2#3{\trivlist
      \item[\hskip \labelsep{\bfseries #2\ #1\ [#3].}]\itshape}
\def\@endtheorem{\endtrivlist}
\makeatother
\def\Z{{\mathbb Z}}

\def\R{{\mathbb R}}
\def\CA{{\cal A}}
\def\CB{{\cal B}}

\def\CM{{\cal M}}
\def\CO{{\cal O}}
\def\CG{{\cal G}}

\def\CI{{\cal I}}
\def\nmo{\hbox{\it n-1}}
\def\dnmo{\hbox{\it 2n-1}}
\def\id{\hbox{id}}
\newcommand\beq{\begin{equation}}
\newcommand\eeq{\end{equation}}
\newcommand{\bea}{\begin{eqnarray}}
\newcommand{\eea}{\end{eqnarray}}
\def\ts{\otimes}

\def\proof{\noindent {\bf Proof:\ }}
\def\endproof{\vrule height 0.5em depth 0.2em width 0.5em}
\newtheorem{lemma}{Lemma}[section]

\newtheorem{proposition}[lemma]{Proposition}
\newtheorem{corollary}[lemma]{Corollary}
\theoremstyle{definition}
\newtheorem{example}[lemma]{Example}
\newtheorem{definition}[lemma]{Definition}
\newtheorem{problem}[lemma]{Problem}
\newtheorem{obse}[lemma]{Observation}
\title{On the $n$-ary algebras, semigroups and their universal covers.}

\author{Andrzej Sitarz\thanks{\mbox{Supported partially by 
bourse post-doctorale du MENRT and KBN grant 2P03B 023 14.}} \\ \ \\
{\small  Laboratoire GCR-CNRS URA 769 - Universit\'e Pierre et Marie Curie,} \\
{\small Tour 22, - 4, Place Jussieu, 75252 Paris Cedex 05, France} \\
{\em and} \\
{\small Institute of Physics, Jagiellonian University,} \\
{\small Reymonta 4, 30059 Krak\'ow, Poland}}

\begin{document}

\maketitle
\begin{abstract}
For any $n$-ary associative algebra we construct 
a $\Z_{n-1}$ graded algebra, which is a universal object
containing the $n$-ary algebra as a subspace of elements 
of degree $1$. Similar construction is carried out for 
semigroups.
\end{abstract}


\section{Introduction}

Recently there has been some interest in the studies of generalized
algebraic structures, in particular, linear spaces, which are equipped
with $n$-ary products. Such objects are natural generalizations 
of algebras. The simplest nontrivial examples, ternary structures, 
have attracted the attention due to some physical motivation 
[MV,VK,K], but has also been a intense subject of studies in 
mathematics (see for instance [Gn,Br] and references therein).

Let us remind here briefly the basic definition.

\begin{definition}
An $n$-ary algebra $\CA$ is a linear space with a linear map 
$m: \CA^{\ts n} \to \CA$. We shall say that $\CA$ is 
an associative\footnote{This notion of associativity is sometimes
denoted as {\em full associativity}.}
$n$-ary algebra if the composition of any $\dnmo$ 
elements is uniquely defined, i.e.:
\begin{equation}
\begin{split}
m ( \id \ts \cdots \ts \id & \ts m_{(i)}  \ts \id \ts \cdots \ts \id) =\\
&=  m ( \id \ts \cdots \ts \id  \ts m_{(j)}  \ts \id \ts \cdots \ts \id)
\end{split}
\end{equation}
for all $1 \leq i,j \leq n-1$, where $m_{(i)}$ denotes that $m$ is on 
the $i$-th position in the tensor product.
\end{definition}

Of course, having a $n$-ary structure one may easily construct 
a $(\dnmo)$-ary one on the top of it and, especially, every ordinary
algebra gives rise to $n$-ary structures (associative or not). 

In this note we shall study the inverted problem: of finding
an algebra, which has a subspace stable under the $n$-ary 
multiplication. We shall demonstrate that all $n$-ary 
associative algebraic structures are of this form and
we shall prove the universality of the constructed object. 

A natural examples of such objects come from 
$\Z_{n-1}$-graded algebras and the subspaces of elements 
of degree $1$ (which are, of course, stable under $n$-multiplication). 
We shall show that the constructed algebra, which covers a given
$n$-ary algebra, is $\nmo$-graded.

Our results open new possibilities for studies of $n$-ary 
objects as well as graded algebras and we shall indicate 
few of them.

\section{Universal $\Z_{n-1}$-graded algebras for $n$-ary algebras.}

Let $\CA$ be an $n$-ary algebra. For a while we do not assume
anything more about $m$.
Let us take $T(\CA)$, the tensor algebra of $\CA$ and the natural
inclusion map $i: \CA \to T(\CA)$. So far, we have embedded $\CA$
(as a linear space) in a $\Z$-graded algebra but only as a linear space.
The following proposition allows us to construct the embedding as 
an $n$-ary homomorphism. 

\begin{proposition}
Let $\CI$ in $T(\CA)$  be a two-sided  ideal generated by the elements 
of the form: $a_1 \otimes \cdots \otimes a_n - (a_1 a_2 \cdots a_n)$.
Let $\CO(\CA) = T(\CA)/\CI$ be the quotient algebra and 
$\xi: T(\CA) \to \CO(\CA)$ the corresponding canonical projection.  

Then, if the $n$-ary algebra $\CA$ is associative, the restriction of 
the projection map $\xi: \CA \to \CO(\CA)$ is an injective 
homomorphism of $n$-ary algebras.
\end{proposition}

\proof Suppose for a while that $\xi$ is injective. Then $\xi(\CA)$
has a natural $n$-ary algebraic structure obtained by taking the
the product in $T(\CA)$ and projecting it again to $\CO(\CA)$. 
{}From the construction of the ideal $\CI$ and projection $\xi$ it
is obvious that the product will be in $\xi(\CA)$ and that $\xi$
is an homomorphism between $\CA$ and $\xi(\CA)$.

Now, suppose that $\xi$ is not injective. Then, there exists an element 
$a \in \CA$, which belongs to the ideal $\CI$. Of course $a$ 
cannot belong to the linear span of generators, therefore it must be of 
the form $a=\sum x \otimes y \otimes z$, where $y$ is in the linear 
span of the generators and $x,z$ are of the form 
$b_1 \otimes \cdots \otimes b_l$, $b_i \in \CA$ (case $l=0$ included).  
Next, we can look at the expression at a fixed degree. Obviously, since 
$a$ is of degree $1$, elements of all other degrees must add up to zero. 
Moreover, we could concentrate our efforts only on elements of degree 
$1$ modulo $\nmo$. 

We shall perform the proof in two steps. First, note that since $a \in \CI$ 
then the part of $\sum x \otimes y \otimes z$ of degree $0$ 
must be of the form $a=\sum_i (a_1^{(i)} \cdots a_n^{(i)})$. Immediately we 
get that in the part of the sum of degree $n$ we must have $\sum_i - (a_1^{(i)} 
\otimes \cdots \otimes a_n^{(i)})$. But the latter do not add up to zero
(in fact we can safely assume that all are linearly independent), so there must
be further components contributing to the sum of all elements of degree $n$. 
These can only come from the expressions of the type $\sum x \otimes y \otimes z$, 
with $y$ being a generator of $\CI$ and $x$,$z$ such that their degrees
add up to $\nmo$. This means that for each $i$ there is a $k$ that $a_k^i$
is again a product of $n$ elements: $a_k^i = (a_k^i)_1 \cdots (a_k^i)_n$. 

We can now repeat the entire procedure step by step going from $n$ to 
$\dnmo$ and further on. At each step either the sum of the elements
of given degree vanishes or we can still go up. However, since our sum 
is finite there is a maximum degree of it and our procedure must stop
at a given moment. 

Then we have the following situation: in each step we have decomposed one
of the elements of $\CA$ as a product of $n$ elements. Such refining goes on
until at a certain degree we arrive at a sum of the type:
$$ \sum a_1 \otimes a_2 \otimes \cdots \otimes (a_k^1 a_k^2 \cdots a_k^n)  
\otimes \cdots \otimes a_p, $$
where $p=m(n-1) +1$ for certain $m$. However, we know that for term of the
above type in the sum the following tensor product:
$$ a_1 \otimes a_2 \otimes \cdots \otimes a_k^1 \otimes a_k^2 \otimes 
\cdots \otimes a_k^n \otimes \cdots \otimes a_p, $$ must be the same. 

Then we can see that from the original expression for $a=\sum_{(i)} 
(a_1^{(i)} \cdots a_n^{(i)})$ 
we have obtained that it could be rewritten as a sum of products of elements of $\CA$ of
length $n$, $\dnmo$ up to $m(n-1) +1$ such that for sums of equal lengths the 
terms differ only in the order in which the product is taken 
and the constant coefficients in front of them add up to $0$. So, if 
our algebra was $n$-ary associative then the order of product does
not matter and the sum vanishes, so $a=0$ and the map $\xi$ is injective. 
\endproof

\proof Clearly, it is sufficient to study the structure of the ideal $\CI$.
Let us note that its generators are all of degree $1$ mod $\nmo$
and therefore $\CI$ becomes naturally a $\Z_{n-1}$-graded subalgebra 
of $T(\CA)$ (we just take the degree of the elements of $\CI$ to be
degree in $T(\CA)$ mod $\nmo$). 

Then, the quotient has again a natural $\Z_{n-1}$-graded structure.
\endproof

We have constructed an algebra, which contains the $n$-ary algebra
as a subspace stable under $n$-multiplication. Since in our construction
we have used the tensor algebra, which carries a natural $\Z$-grading,
we might ask whether $\CO{\CA}$ is a graded algebra:

\begin{corollary}
The algebra $\CO(\CA)$ is a $\Z_{n-1}$ graded associative algebra
and the subspace of elements of grade $1$ is $n$-ary isomorphic 
with $\CA$.
\end{corollary}

\begin{corollary}
Suppose now that $\CB$, $\CB  \subset \CA$ is a $n$-ary associative
subalgebra of $\CA$. Then $\CO(\CB)$ is  a $\nmo$ graded subalgebra 
of $\CO(\CA)$.
\end{corollary}

The proof of the corollary is simple: first one observes that $T(\CB) \subset T(\CA)$
and, moreover $\CI_\CB \subset \CI_\CA$. Then the appropriate relation 
for the inclusion of $(\nmo)$-graded algebras follows.

Now, we shall state the main proposition, which concerns the universality property of
the constructed object.

\begin{proposition}
Let $M$ be  a $\Z_{n-1}$-graded associative algebra, and let $M_1$ denote the space of all its elements of degree $1$, which is then an $n$-ary 
associative algebra.  If  there exists an a $n$-ary homomorphism 
\mbox{$\rho: \CA \to M_1 \subset M$} then there exist a unique 
homomorphism 
of $\Z_{n-1}$-graded algebras $\tilde{\rho}:  \CO(\CA) \to M$ such that 
the following diagram is commutative:

$$\xymatrix{  \CA  \ar[d]_i \ar[r]^{\rho}  &  M \\ 
 \CO(\CA)   \ar[ur]_{\tilde{\rho}} } $$
\end{proposition}

\proof First, note that the map $\rho$ is linear. Let us extend it in
a natural way to $T(\CA)$ by taking:
$$T(\rho)( a_1 \otimes \cdots a_k) = \rho(a_1) \cdots \rho(a_k).$$

Now the only thing is to check that $T(\rho)$ vanishes on the ideal 
$\CI$, but due to the fact that $\rho$ is an $n$-ary homomorphism 
this is obvious. Similarly, by construction it is also clear that $\hat{\rho}$ 
preserves the grading.  \endproof

What we have shown are two important facts: first, that every $n$-ary
associative algebra could be identified with the space of elements
of degree $1$ of some $\Z_{n-1}$-graded algebra, second that the latter
is universal, i.e. for every $\Z_{n-1}$-graded algebra in which our $n$-ary
structure is embedded there exists a homomorphism between them. 

Therefore the above lemma solves the problem posed in the introduction:
not only we know that every $n$-ary algebra can be embedded in 
a {\em normal} (binary) algebra but we know how to find such objects.

\begin{corollary}
Let $\CA$, $M$ and $\rho: \CA \to M$ be as defined in the universality 
proposition. Let us call $\hat{M}$ the subalgebra of $M$ generated by 
the image of $\rho$. Then $\hat{M}$ is isomorphic (as a $\Z_{n-1}$-graded 
algebra) to $\CO(\CA)/I$ for some ideal $I \subset \CO(\CA)$.
\end{corollary}

\begin{corollary}
Let $\CA$ and $\CB$ be $n$-ary associative algebras and 
\mbox{$\phi: \CA \to \CB$} a homomorphism. Then there exists a unique 
homomorphism  of  \mbox{$\Z_{\nmo}$-graded} algebras 
$\CO(\phi): \CO(\CA) \to \CO(\CB)$ such that the following diagram 
is commutative:

$$\xymatrix{ \CA \ar[r]^{\phi}  \ar[d]_\xi & \CB \ar[d]_\xi \\ 
\CO(\CA) \ar[r]^{\CO(\phi)} & \CO(\CB)  }$$
\end{corollary}

where, $\xi$ denotes the embedding of $\CA,\CB$ into
$\CO(\CA),\CO(\CB)$, respectively. Again, the proof of both 
corollaries is a simple consequence of the universality lemma.

\section{Examples.}

In the previous section we have learned how to construct 
a $\Z_{n-1}$-graded universal envelope algebra of the $n$-ary algebra $\CA$. 
A simple example of ternary ($n=3$) structure would 
illustrate the problem. 

\begin{example}
Let us take a ternary algebra $\CA$ of the odd degree elements of the exterior 
tensor algebra of $\R^n$ (with the product obtained from the exterior product). 
Then its $\Z_2$-graded universal enveloping algebra is equal 
$\CA \oplus \CA_1 \otimes \CA$, where $\CA_1$ denotes the linear 
span of the generators of $\CA$. 
\end{example}

The constructed universal embedding is rather big - note that in the above example 
we do not recover the exterior algebra over $\R^n$ we have started with but a space 
much larger. The following observation suggests the way to restrict it:

\begin{obse}
Let $\CA$  be an associative $n$-ary algebra, $\CO(\CA)$ its $\Z_{n-1}$-graded
envelopping algebra $i: \CA \to \CO(\CA)$ and $I \subset \CO(\CA)$ be a graded 
ideal which does not intersect the image of $\CA$. Then  $\CA$ can be embedded 
into $\CO(\CA)/I$
\end{obse}

Of course, out of the mentioned class of ideals with this property we can always find 
a maximal ideal, then we shall obtain the {\em minimal} $\Z_{n-1}$-graded 
algebras which contain $\CA$. 

The maximal ideal, however, could be too big and in fact we might not recover
the original algebra we have started with. Let us illustrate it with an example.

\begin{example}
For the above simple ternary example we might  choose the ideal generated by
the symmetric part of the tensor product $\CA_1 \otimes \CA_1$. Note, that 
for an even $n$ this is {\bf not} a maximal ideal: we can add to it an element, 
which is of the form $e^{1} \otimes e^{2}  e^{3} \cdots e^{n}$, as it is annihilated 
by the action of the entire algebra from both sides and clearly does not belong 
to our previously chosen ideal.

For our original choice we recover as the quotient the original $\Z_2$-graded 
exterior algebra and in the case of  the extended (maximal) ideal we would have
the product of $n$ generators $e^{1} \otimes e^{2}  e^{3} \cdots e^{n}$ vanishing.
\end{example}

\section{Semigroups}

A vast class of algebras come from groups, constructed as groups algebras
and, similarly, $\Z_n$-graded algebras can originate from $\Z_n$-graded 
semigroups. Let us propose a definition.

\begin{definition}
Let $\CG$ be a set with a map $m: \CG^n \to \CG$, then we shall call
it an $n$-semigroup. 
\end{definition}

Of course, there is a difficulty in extending this notion, for instance, 
by introducing a unit element in such a way that it does not reduce
the $n$-product to the binary one. 

Similarly as in the case of algebras some nontrivial examples of 
$n$-semigroups come from $Z_{n-1}$-graded groups (groups with a 
given homomorphism on $\Z_{n-1}$) by considering the inverse 
image of $1$.

Now, we can pose the following question:  can every $n$-ary associative 
semigroup be embedded into a semigroup and, what are the conditions 
that allow it to be embedded in a group?

\begin{proposition}
If $\CM$ is a associative $n$-ary semigroup then there exists a $\Z_{n-1}$ graded 
semigroup $\CO(\CM)$ such that $i: \CM \to \CO(\CM)$ is an isomorphism
between $\CM$ and the subspace of elements of degree $1$. 

Moreover the construction is universal in the following sense: for every $\Z_n$-graded 
semigroup $N$  and a homomorphism $\rho: \CM \to N$ there 
exist a homomorphism $\hat{\rho}$ such that the following diagram is
commutative:
$$\xymatrix{  \CM  \ar[d]_i \ar[r]^{\rho}  &  N \\ 
\CO(\CM)   \ar[ur]_{\tilde{\rho}} } $$
\end{proposition}

\proof Although the subject of the lemma is similar as in the case
of algebras the problem is slightly more complicated. Again we begin 
by constructing the direct sum $T(\CM)$ of $\CM^k$ (cartesian 
product of $k$ copies of $\CM$), for $k>0$. We introduce the product 
in this space as the standard cartesian product. 
Now we may introduce a relation in $T(\CM)$ - we say that two elements
are related if the product of all their components are equal to each other:
$$  g_1 \times \cdots \times g_k \sim h_1 \times \cdots \times h_l 
\hbox{\ \ iff \ \ } g_1 g_2 \cdots g_k = h_1 h_2 \cdots h_k. $$

This relation is clearly an equivalence relation (note that it makes sense
only for \mbox{$k,l=r_{k,l}(n-1)+1$} and is well-defined due to associativity).
We extend it that it agrees with the product, i.e. for any 
$x,y,w,z \in T(\CM)$ and $x \sim y$ we postulate \mbox{$z \times x 
\times w \sim z \times y \times w$}. 

Therefore we might introduce the quotient of $T(\CM)$ by the relation 
$\sim$ and transport the product to the quotient 
$\CO(\CM) = T(\CM)/\sim$.

The map $i : \CM \to \CO(\CM)$ is clearly injective (from the
definition no elements of $\CM$ could be related among themselves)
and preserving the $n$-ary product. What remains to be checked is
the universality property and $\Z_{n-1}$-grading. This, however,
almost exactly copies the idea of the proof for $n$-ary algebras.
\endproof

Let us notice that we were able to proof the correspondence between
$n$-semigroups and semigroups and the construction cannot tell us
whether the universal object (or the $n$-semigroups we started from)
is related with a monoid or a group. We shall come to this problem
later.

\section{Properties of $n$-ary algebras and $n$-semigroups}

One of the consequences of the proven correspondence between
$n$-ary objects and $\Z_{n-1}$-graded algebras is the possibility
to translate several constructions proposed for $n$-ary algebras
(semigroups). 

We shall indicate here two problems: one of Hochschild homology,
and another related with ternary semigroups.

\subsection{Hochschild Homology of $n$-ary algebras}

A generalization of Hochschild homology  has been proposed in 
some papers (see [MV,Gn]), some modification will also be 
discussed in [Si]. Let us remind the definition [Gn]:

\begin{definition}
Let $\CA$ be an $n$-ary algebra, for an even $n$. Let $C^k$ denote
$\CA^{k(n-1)+1}$. Consider a map $\delta_i^k : C^k \to C^{k-1}$
defined as follows:
\begin{equation}
\begin{split}
\delta_i (a_0 \ts \cdots \ts a_i \ts &\cdots \ts a_{i+n-1}\ts \cdots
\ts a_{k(n-1)} ) = \\
&= (a_0 \ts \cdots \ts ( a_i a_{i+1}\cdots a_{i+n-1}) \ts \cdots
\ts a_{k(n-1)}.
\end{split}
\end{equation}
Then $d_k = \sum_i (-1)^i \delta_i^k: C^k \to C^{k-1}$ is a linear
operator which satisfies $d_k d_{k+1} = 0$. (For details and proof 
of this statement see [Gn]).
\end{definition}

Now, let us formulate the problem:

\begin{problem}
How does the above defined homology of the chain complex 
$(C^k(\CA), d)$ relate to the usual Hochschild homology of 
the universal cover of $\CA$? 

Is there a chain complex and homology of $\CO(\CA)$
(also a generalized homology with a boundary satisfying $d^N=0$
for some $N$, see [DV] for details) such that Hochschild homology
of $\CA$ can be expressed in terms of this homology?
\end{problem}

\subsection{Ternary groups}

By looking at the elements of degree $1$ of a $\Z_{n-1}$-graded
group it is obvious that they do not form a subgroup, nevertheless 
still some specific group operations do exist, for instance, for every
element one can find $n-2$ elements such that their product gives
the unit of the groups. 

A very special situation occurs when we have $n=3$ - then the 
operation is unique since the map $g \to g^{-1}$ does not change
the degree of an element.

For this reason when we take this very simple case $n=3$ we might
be able to propose a definition of a ternary group:

\begin{definition}
A ternary groups $\CG$ is a set, with a ternary associative map 
$\CG^3 \to \CG$ and the inverse map $\CG \ni g \to g^{-1} \in 
\CG$ such that:
\begin{eqnarray*}
&& \forall g,h \in \CG \;\; gg^{-1}h = h = hgg^{-1}
\end{eqnarray*}
\end{definition}

Notice that for an ordinary group, the above statement is equivalent to
the existence of a unit and an inverse. Here, however, we cannot draw
the same conclusion. We begin with a very easy lemma.

\begin{lemma}
For a ternary groups the map $h \to g^{-1} h g$ is an injective
homomorphism.
\end{lemma}

\proof It is clear that it is a morphism, the injectivity follows from the
existence of its inverse $h \to g h g^{-1}$. \endproof

Notice that unlike in the classical case this does not proof that it is
surjective. Of course in the finite dimensional case (i.e. when cardinality
of the set $G$ is finite). Now we can state the problem:

\begin{problem}
Are there non-trivial ternary groups (in the above sense), which cannot
be embedded in usual ($\Z_2$-graded) groups?
\end{problem}   

\section{Conclusions}

In this paper we have shown that the problem of at least associative $n$-ary
algebras can be reduced to the {\nmo}-graded usual algebras. This has many
implications: first of all, one can translate the result of studies of $n$-ary 
objects to the language of graded algebras, as we have suggested here in the
case of Hochschild homology. Since our results concern only associative 
structures, one may try investigate whether analogous relations are present
in the arbitrary case, of Lebniz, or Lie-type structures.

Another application comes from the opposite direction. For $n$-ary algebras
one can have several notions of commutativity or generalizations of 
anticommutativity. An example of that is $j$-commutativity of ternary algebras,
for $j$ being a cubic root of unity, for any $a,b,c \in \CA$ we impose:
$$ abc = j bca = j^2 cab.$$

This relation can be easily translated for relations between elements of the
universal (or any other) $\Z_2$-graded covering of $\CA$. This provides us 
with a new class of algebraic objects, with interesting commutation relations.

{\bf Acknowledgement}: The author would like to thank M.Dubois-Violette
and R.Kerner for discussions and helpful remarks, and IHES for hospitality.

\end{document}